\documentclass{amsart}

\usepackage{amssymb}
\usepackage[margin=1in]{geometry}
\usepackage{graphicx}
\usepackage[dvipsnames]{xcolor}
\usepackage{xurl}
\usepackage[
  hypertexnames=false,
  colorlinks=true,
  linkcolor=Blue,
  citecolor=Blue,
  urlcolor=Blue
]{hyperref}
\usepackage[nameinlink,capitalise,noabbrev]{cleveref}
\usepackage[inline,shortlabels]{enumitem}
\usepackage{aliascnt}
\usepackage[
  backend=biber,
  style=alphabetic
]{biblatex}
\usepackage{booktabs}
\usepackage{tabularx}
\usepackage{array}
\usepackage{float}
\usepackage{comment}

\theoremstyle{plain}
\newtheorem{theorem}{Theorem}[section]
\newtheorem{example}[theorem]{Example}

\newcommand{\newsharedtheorem}[2]{%
  \newaliascnt{#1}{theorem}%
  \newtheorem{#1}[#1]{#2}%
  \aliascntresetthe{#1}%
  \crefname{#1}{\MakeLowercase{#2}}{\MakeLowercase{#2}s}%
  \Crefname{#1}{#2}{#2s}%
}

\newsharedtheorem{proposition}{Proposition}
\newsharedtheorem{corollary}{Corollary}
\newsharedtheorem{lemma}{Lemma}
\newsharedtheorem{claim}{Claim}

\theoremstyle{definition}
\newsharedtheorem{definition}{Definition}

\theoremstyle{remark}
\newsharedtheorem{remark}{Remark}

\crefname{theorem}{theorem}{theorems}
\Crefname{theorem}{Theorem}{Theorems}

\DeclareMathOperator{\Aut}{Aut}
\DeclareMathOperator{\rk}{rk}

  \DeclareMathOperator{\lcp}{lcp} 

\newcommand{\bbN}{\mathbb{N}}
\newcommand{\bbP}{\mathbb{P}}
\newcommand{\bbF}{\mathbb{F}}

\newcommand{\cF}{\mathcal{F}}

\newcommand{\cG}{\mathcal{G}}

\newcommand{\Range}{{\operatorname{Range}}}
\title{Random Non-Commutative Error-Correcting Codes}
\date{\today}

\author[M.\ Best]{Michael Best}
\address{Michael Best\hfill\break
	Weizmann institute of Science\hfill\break
	Rehovot, Israel.}
\email{michael.best@weizmann.ac.il}

\author[A.\ Lubotzky]{Alexander Lubotzky}
\address{Alexander Lubotzky\hfill\break
	Weizmann institute of Science\hfill\break
	Rehovot, Israel.}
\email{alex.lubotzky@mail.huji.ac.il}

\begin{document}
\begin{abstract}
Chapman, Dinur, and Lubotzky recently introduced the theory of
non-commutative error-correcting codes \cite{ChapmanDinurLubotzky},
generalizing classical coding theory to arbitrary groups via $\delta$-robust
subsets used to detect proper subgroups. Building
on this foundation, this paper systematically studies an analogue of random
codes for arbitrary groups, exploring the extent to which the classical
Gilbert--Varshamov (GV) bound holds in this non-commutative setting. While
random universal codes of linear size in free groups are shown to be
asymptotically bad, we establish a complete characterization for general
families of finite groups. Specifically, a family of finite groups admits
good random codes (forming a GV family) if and only if it exhibits uniform
polynomial maximal-subgroup growth. Furthermore, we provide explicit
quantitative bounds on the rate and optimal distance range based on intrinsic
parameters, proving that families of finite solvable groups always form GV
families and determining the exact probabilistic rate for direct powers of a
fixed finite simple group.
\end{abstract}
\maketitle

\section{Introduction}

Let $p$ be a prime and let $\bbF_p^n$ be the $n$-dimensional vector space
over the finite field $\bbF_p$. An \emph{$[n,k,d]_p$ linear code} is a
linear subspace $\mathcal C$ of $\bbF_p^n$ of dimension $k$ and distance
$d$; that is, the Hamming weight of every nonzero vector
$v\in\mathcal C$ is at least $d$. A family of $[n_i,k_i,d_i]_p$ codes
$\mathcal C_i\leq\bbF_p^{n_i}$, where $n_i\longrightarrow\infty$, is called
\emph{(asymptotically) good} if there exist $\delta,\rho>0$ such that
$d_i/n_i\geq\delta$ and $k_i/n_i\geq\rho$. We call
$\rho(\mathcal C_i):=k_i/n_i$ the \emph{rate} of the code and
$\delta(\mathcal C_i):=d_i/n_i$ its \emph{normalized distance}. The basic
result of Shannon asserts:

\begin{theorem}[Shannon]
\label{thm:good-codes}
\cite[Theorem~11]{Shannon1948}.
There exist good codes.
\end{theorem}

The second basic result of the theory of error-correcting codes is the
Gilbert--Varshamov inequality.

\begin{theorem}[Gilbert--Varshamov bound]
\label{thm:GV}
\cite[Theorem~1, p.~507]{Gilbert1952} and
\cite[condition~(3), pp.~739--741]{Varshamov1957}.
Let $0<\delta<1-1/p$ and $0<\rho<1-H_p(\delta)$. Then a random linear
code in $\bbF_p^n$ of rate $\rho$ has normalized distance at least
$\delta$ with high probability.
\end{theorem}

Here, for a real number $s>1$, we denote the \emph{$s$-ary entropy} by
$H_s(\delta)$ for $0<\delta<1$, where
\begin{equation}
\label{eq:s-ary-entropy}
H_s(\delta)
=\delta\log_s(s-1)-\delta\log_s\delta
 -(1-\delta)\log_s(1-\delta).
\end{equation}
By a \emph{random code} $\mathcal C$ of rate $\rho$, we mean that
$k=\lfloor\rho n\rfloor$ vectors in $\bbF_p^n$ are chosen independently
and uniformly at random and that $\mathcal C$ is the subspace they span.
The condition $\delta<1-1/p$ is necessary because of the Plotkin bound.

\begin{theorem}[Plotkin bound]
\label{thm:plotkin}
\cite[Theorems~1--2, pp.~446--447]{Plotkin1960} and
\cite[Theorem~4.4.1]{GuruswamiRudraSudanECT}. If
$\mathcal C_i\leq\bbF_p^{n_i}$, $n_i\longrightarrow\infty$, and
$\delta(\mathcal C_i)\geq1-1/p$, then
$\rho(\mathcal C_i)\longrightarrow0$.
\end{theorem}

In \cite{ChapmanDinurLubotzky}, Chapman, Dinur, and the second author
established a far-reaching extension of \Cref{thm:good-codes}. To state it,
we first need the following definition.

\begin{definition}
\label{def:delta-robust}
Let $G$ be any finitely generated group, let $0<\delta<1$, and let
$\mathcal A$ be a finite tuple in $G$, viewed as a multiset. We say that
$\mathcal A$ is \emph{$\delta$-robust} if every submultiset containing more
than $(1-\delta)|\mathcal A|$ of its entries generates $G$. All cardinalities
and intersections involving sampled tuples are understood to count
multiplicity.
\end{definition}

For a finite tuple $\mathcal A$ in $G$, its \emph{detection probability} is
\[
     \delta(\mathcal A;G)
:=
\inf
\left\{
1-\frac{|\mathcal A\cap H|}{|\mathcal A|}: H\lneqq G
\right\}.
\]
Thus $\mathcal A$ is $\delta$-robust if and only if
$\delta(\mathcal A;G)\geq\delta$. We write $\rk(G)$ for the minimal number
of generators of $G$.

Chapman, Dinur, and the second author established the following
non-commutative analogue of Shannon's theorem.

\begin{theorem}
\label{thm:good-univ-codes}
\cite[Theorem~1.8]{ChapmanDinurLubotzky}.
There are absolute constants $\delta>0$ and $C>1$ such that every
finitely generated group $G$ has a $\delta$-robust subset of size at most
$C\rk(G)$.
\end{theorem}

\Cref{thm:good-codes} is a special case of
\Cref{thm:good-univ-codes}. Indeed, apply
\Cref{thm:good-univ-codes} to $G=\bbF_p^k$ to obtain a robust subset
$\mathcal A=\{\alpha_1,\ldots,\alpha_n\}\subseteq\bbF_p^k$. Let
$\mathcal C$ be the image of the linear map
\[
T:\bbF_p^k\longrightarrow\bbF_p^n,
\qquad
T(\beta)
=
(\langle\beta,\alpha_1\rangle,\ldots,
 \langle\beta,\alpha_n\rangle).
\]
Here, for $\alpha=(a_1,\ldots,a_k)$ and
$\beta=(b_1,\ldots,b_k)$, we put
$\langle\beta,\alpha\rangle=\sum_{j=1}^k a_jb_j$. One easily checks that
$\mathcal C=T(\bbF_p^k)\leq\bbF_p^n$ has dimension $k$ and normalized
distance at least $\delta$. Since $n\leq Ck$, its rate is at least $1/C$.

Every $k$-generated group is a quotient of the free group $F_k$ of rank
$k$, and
robustness is preserved under quotients. It therefore suffices to construct
robust subsets in $F_k$; such subsets are called \emph{universal codes}.
The goal of this paper is to study random non-commutative codes
systematically, and in particular to determine the extent to which an
analogue of \Cref{thm:GV} holds for general groups.

Recall that a random linear code may be obtained from a random
$k\times n$ generator matrix, considered in the binary case in
\cite[Section~II-C]{BargForney2002}. Viewing the matrix by columns gives a
random $n$-tuple in $\bbF_p^k$. This suggests our general model: we ask when
a random tuple in a group $G$ is robust.

To formalize this model, let
$\cG=(G_i)_{i\geq1}$ be a
family of nontrivial finitely generated groups, where each $G_i$ is equipped
with a probability measure $\mathbb P_i$. For $0<\rho<1$, put
$
N_i(\rho):=\left\lceil\rk(G_i)/\rho\right\rceil
$
and let $A_i(\rho)$ be a tuple of $N_i(\rho)$ independent
$\mathbb P_i$-distributed elements of $G_i$. We use the same notation
$\mathbb P_i$ for the corresponding product measure on such tuples.

\begin{definition}[$(\delta,\rho)$-GV family]
\label{def:probabilistic-GV}
The family $\cG$ is a \emph{$(\delta,\rho)$-Gilbert--Varshamov family},
or a \emph{$(\delta,\rho)$-GV family}, if
\[
\bbP_{\cG}(\delta,\rho)
:=
\liminf_{i\to\infty}
\mathbb P_i\!\left[\delta(A_i(\rho);G_i)\geq\delta\right]
>0.
\]
We call $\cG$ a \emph{GV family} if this holds for some
$0<\delta,\rho<1$.
We call $\cG$ \emph{strong GV} if, for every $0<\delta<1$, there is some
$0<\rho<1$ for which $\cG$ is $(\delta,\rho)$-GV.
\label{def:strong-GV}
\end{definition}

In other words, a GV family admits random robust tuples whose size is
proportional to the generating rank, with a success probability bounded away
from zero along the family. Being strong GV means that this remains possible
for every target normalized distance below $1$, although the rate may depend
on the target distance.

There are two rather different settings. For free groups one must first
choose a probability model, and the natural models considered here give a
negative answer. For finite groups there is a canonical uniform measure,
and the answer is governed by maximal-subgroup growth. These are
the two main results of the paper.

\subsection{Random universal codes}
\label{subsec:random-unic-codes-intro}
The first result concerns a question suggested in
\cite[Problem~7.1]{ChapmanDinurLubotzky}. For
$F_k=\langle x_1,\ldots,x_k\rangle$ and $S=\{i_1<\cdots<i_r\}\subseteq[k]$,
put $x_S=x_{i_1}\cdots x_{i_r}$ and
$
\mathcal H_k:=\{x_S:S\subseteq[k]\}.
$
As shown in \cite[Theorem~3.3 and its proof]{ChapmanDinurLubotzky}, this
universal Hadamard code has $2^k$ elements and detection probability at
least $1/2$. Thus its normalized distance is bounded away from zero,
but its rate $k/2^k$ tends to zero. In the classical abelian Hadamard
setting, taking a random linear-size sample produces a good code. It was
natural to ask whether the same procedure works for $\mathcal H_k$.

\begin{theorem}
\label{thm:random-hadamard-not-robust-intro}
Let $\mathbb P_k$ be the uniform probability measure on $\mathcal H_k$.
For $C>1$, let $A_k(C)$ be a tuple of $\lceil Ck\rceil$ independent
$\mathbb P_k$-distributed elements of $F_k$. Then
\[
\bbP_k\!\left[\delta(A_k(C);F_k)=0\right]
\xrightarrow[k\to\infty]{}1.
\]
In particular, the family $(F_k,\mathbb P_k)_{k\geq1}$ is not GV.
\end{theorem}

Thus a random linear-size sample from the universal Hadamard code is not a
good universal code. This gives a negative answer to the question above,
despite the positive normalized distance of the ambient Hadamard code.
We obtain the same negative result for a linear-size random sample from
word-metric balls;
see
\Cref{sec:random-univ-codes}.

\subsection{Random codes in finite groups}

For the remainder of this subsection, and throughout
\Cref{sec:finite-groups}, the groups $G_i$ are finite and, unless specified
otherwise, $\mathbb P_i$ is the uniform probability measure on $G_i$.

To state our results for GV families of finite groups, write $q(G)$
for the smallest index of a proper subgroup of $G$, and put
$q_{\cG}:=\liminf_{i\to\infty}q(G_i)$.

We say that $\cG$ has \emph{uniform polynomial maximal-subgroup growth},
or \emph{uniform PMSG}, if there is a constant $A\geq0$ such that,
for all sufficiently large $i$ and every $s\geq2$, the number of
maximal subgroups of index $s$ in $G_i$ is at most
$s^{A\rk(G_i)}$. Thus the polynomial exponent is controlled uniformly
by the generating rank.

\begin{theorem}[GV characterization]
\label{thm:GV-characterization-intro}
Let $\cG=(G_i)_{i\geq1}$ be a family of nontrivial finite groups, and adopt
the convention that $1/q_{\cG}=0$ if $q_{\cG}=\infty$. Then the following
are equivalent:
\begin{enumerate}
    \item The family $\cG$ is $(\delta,\rho)$-GV for some
    $\delta>0$ and $0<\rho<1$.
    \item The family $\cG$ has uniform PMSG.
    \item For every $0<\delta<1-1/q_{\cG}$, there exists
    $0<\rho<1$ such that $\cG$ is $(\delta,\rho)$-GV.
\end{enumerate}
\end{theorem}

This gives a complete group-theoretic characterization, rather than only
a sufficient construction. The connection between maximal-subgroup growth
and random generation was first observed by Lubotzky
\cite{Lubotzky1995} and studied systematically by Mann
\cite{Mann1996}, who developed the theory of positively finitely generated
profinite groups. Mann and Shalev \cite{MannShalev1996} proved that positive
finite generation is equivalent to polynomial maximal-subgroup growth for
profinite groups. Lucchini and Moscatiello \cite{LucchiniMoscatiello}
provide the quantitative input used here.

To state our quantitative results, we record the rate and range of a
family and introduce the technical parameters used here. We recall their
precise maximal-subgroup definitions in \Cref{sec:finite-groups}.

\begin{definition}[GV rate, range, and family parameters]
\label{def:GV-rate}
\label{def:family-parameters}
For $0<\delta<1$, define
\[
R_{\cG}(\delta)
:=
\sup\{\rho\in(0,1):\cG\text{ is }(\delta,\rho)\text{-GV}\},
\qquad
\Range(\cG)
:=
\sup\{\delta\in(0,1):R_{\cG}(\delta)>0\}.
\]
The supremum of the empty set is understood to be zero. For a nontrivial
finite group $\Gamma$, define, for each $s\geq2$,
\[
m_s(\Gamma):=\#\{M<_{\max}\Gamma:[\Gamma:M]=s\},
\text{ and set }
\mathcal M(\Gamma):=
\max_{\substack{s\geq2\\m_s(\Gamma)>0}}
\frac{\log m_s(\Gamma)}{\log s}.
\]
Put
$r_i:=\rk(G_i)$ and define
\[
\mu_{\cG}:=\max\left\{1,\limsup_{i\to\infty}
\frac{\mathcal M(G_i)}{r_i}\right\},\qquad
\rk_{\cG}:=\liminf_{i\to\infty}r_i\in\bbN\cup\{\infty\}.
\]
The parameter $\mu_{\cG}$ is the optimal exponent in the uniform PMSG
bound: it is the infimum of admissible exponents $A\geq1$ for which
$m_s(G_i)\leq s^{A r_i}$ for every $s\geq2$ and all sufficiently large
$i$, with $\mu_{\cG}=\infty$ if no such exponent exists. In other words,
$\mu_{\cG}$ measures the optimal polynomial growth rate of maximal
subgroups relative to the generating rank. Thus $\cG$ has uniform PMSG
if and only if $\mu_{\cG}<\infty$.
\end{definition}

For a fixed prime $p$, \Cref{thm:GV} says that
$\cG_p:=(\bbF_p^k)_{k\geq1}$ is a GV family, while \Cref{thm:plotkin} shows that
it is not strong GV. In fact,
\[
    R_{\mathcal{G}_p}(\delta)
    =
    \begin{cases}
    1-H_p(\delta),
        & 0< \delta<1-\dfrac{1}{p}, \\[4pt]
    0,
        & 1-\dfrac{1}{p}\leq\delta<1,
    \end{cases}
    \qquad
    \Range(\mathcal{G}_p)
    =
    1-\frac{1}{p}.
\]
This classical family is the prototype for the
finite-group theory developed below.

We also characterize the strong GV property; see
\Cref{prop:strong-GV}. Beyond these qualitative characterizations, we obtain
quantitative bounds on the GV rate and determine the optimal distance range
in the relevant regimes. The formal parameters and precise statements are
given in \Cref{sec:finite-groups}, and \Cref{tab:GV-range-summary} summarizes
the resulting characterization, distance ranges, and rate bounds.

\begin{table}[H]
\centering
\small
\renewcommand{\arraystretch}{1.35}
\begin{tabularx}{\textwidth}{
    c c c
    >{\raggedright\arraybackslash}p{2.4cm}
    c
    >{\raggedright\arraybackslash}X
}
\toprule
$\mu_{\cG}$ & $q_{\cG}$ & $\rk_{\cG}$
& GV/strong GV
& $\Range(\cG)$
& $R_{\cG}(\delta)$ \\
\midrule

$\infty$
& arbitrary
& arbitrary
& $\cG$ is not GV.
& $0$
& For every $0<\delta<1$, $R_{\cG}(\delta)=0$. \\

\cmidrule(lr){1-6}

$<\infty$
& $\infty$
& $\infty$
& $\cG$ is strong GV.
& $1$
& For every $0<\delta<1$,
  $R_{\cG}(\delta)=\dfrac{1-\delta}{\mu_{\cG}}$. \\

\cmidrule(lr){1-6}

$<\infty$
& $\infty$
& $<\infty$
& $\cG$ is strong GV.
& $1$
& For every $0<\delta<1$,
  $R_{\cG}(\delta)
  \geq
  \dfrac{1-\delta}{\mu_{\cG}+1/\rk_{\cG}}$. \\

\cmidrule(lr){1-6}

$<\infty$
& $<\infty$
& $\infty$
& $\cG$ is GV but not strong GV.
& $1-\dfrac1{q_{\cG}}$
& For $0<\delta<1-1/q_{\cG}$,
  $\dfrac{1-H_{q_{\cG}}(\delta)}{\mu_{\cG}}
  \leq R_{\cG}(\delta)
  \leq\dfrac{1-\delta}{\mu_{\cG}}$.
  For $1-1/q_{\cG}<\delta<1$,
  $R_{\cG}(\delta)=0$. \\

\cmidrule(lr){1-6}

$<\infty$
& $<\infty$
& $<\infty$
& $\cG$ is GV.
& $\geq1-\dfrac1{q_{\cG}}$
& For every $0<\delta<1-1/q_{\cG}$,
  $R_{\cG}(\delta)
  \geq
  \dfrac{1-H_{q_{\cG}}(\delta)}
  {\mu_{\cG}+2/\rk_{\cG}}$. \\

\bottomrule
\end{tabularx}
\caption{The GV characterization, distance range, and rate bounds in
terms of the family parameters.}
\label{tab:GV-range-summary}
\end{table}

Two applications illustrate the characterization. First, every family of
nontrivial finite solvable groups has uniform polynomial maximal-subgroup
growth and is therefore a GV family. Using
\cite{LucchiniMenegazzoMorigi2006}, we also obtain explicit bounds on the
GV rate; See \Cref{sec:examples} for details.

Second, the threshold can be computed exactly for
direct powers of a fixed finite simple group $H$.
This example includes the classical case. Indeed, if $H$ is abelian, then
$H\cong C_p\cong(\bbF_p,+)$ for some prime $p$, so $H^m\cong\bbF_p^m$ and
random tuples in $H^m$ correspond to random generator matrices for
$p$-ary linear codes. If $H$ is nonabelian, then
$
\rk(H^m)=\log_{|H|}m+O(1)
$
as $m\to\infty$; see \cite[Proposition~3.4]{Thevenaz1997} and
\cite[Lemma~2]{Wiegold1978}. Moreover, the smallest index of a proper
subgroup of $H$ is less than $\sqrt{|H|}$. Despite these differences,
the probabilistic GV rate in both cases is given by the same formula.

\begin{theorem}[Direct powers of a finite simple group]
\label{thm:simple-power-probabilistic-profile-intro}
Let $H$ be a fixed finite simple group and let
$\cG_H=(H^m)_{m\geq1}$. For every
$0<\delta<1$,
\[
R_{\cG_H}(\delta)
    =
    \begin{cases}
    \log_{|H|}\!q(H)\bigl(1-H_{q(H)}(\delta)\bigr),
        & 0<\delta<1-\dfrac{1}{q(H)},\\[4pt]
    0,
        & 1-\dfrac{1}{q(H)}\leq\delta<1.
    \end{cases}
    \qquad
    \Range(\cG_H)=1-\frac{1}{q(H)}.
\]
For every $0<\delta<1-1/q(H)$ and every $0<\rho<1$,
\[
\bbP_{\cG_H}(\delta,\rho)
=
\begin{cases}
1,&\rho\leq R_{\cG_H}(\delta),\\
0,&\rho>R_{\cG_H}(\delta).
\end{cases}
\]
\end{theorem}
Thus direct powers of every finite simple group satisfy the same
entropy-based threshold formula as in the abelian case.

The paper is organized as follows. \Cref{sec:preliminaries} collects the
entropy estimates used throughout. In \Cref{sec:random-univ-codes} we
prove the negative results for the universal Hadamard code and word-metric
balls in free groups. In \Cref{sec:finite-groups} we prove the finite-group
characterization and its quantitative refinements. Finally,
\Cref{sec:examples} treats direct powers, finite solvable groups, and a
family that fails the characterization.
\section*{Acknowledgments}
The authors thank Michael Chapman for his useful comments on this work.
The second author was partially supported by the European Research
Council (ERC) under the European Union's Horizon 2020 research and
innovation programme (grant agreement No.~882751, TeStability).

The authors acknowledge the use of OpenAI's ChatGPT to assist with reviewing and refining proofs, editing and proofreading this work, and identifying relevant references. The authors take full responsibility for the mathematical content.

\section{Preliminaries}
\label{sec:preliminaries}

\begin{lemma}[Binomial entropy tail]
\label{lem:binomial-entropy-tail}
Let $s\geq2$, let $0<\delta\leq1-1/s$, and, for $N\geq1$, put
$
t_N:=\left\lfloor(1-\delta)N\right\rfloor+1.
$
If
$
X_N\sim\operatorname{Bin}\left(N,\frac1s\right),
$
then
\begin{equation}
\Pr(X_N\geq t_N)
=
\frac{V_s(N,N-t_N)}{s^N},
\label{eq:binomial-hamming-identity}
\end{equation}
where $V_s(N,r):=\sum_{k=0}^{r}\binom Nk(s-1)^k$.
So,
\begin{equation}
\Pr(X_N\geq t_N)
\leq
s^{-N(1-H_s(\delta))}.
\label{eq:binomial-entropy-upper}
\end{equation}
For fixed $s$ and $\delta$, the Hamming-ball volume estimate also gives
\begin{equation}
\Pr(X_N\geq t_N)
=s^{-N(1-H_s(\delta))+o(N)}.
\label{eq:binomial-entropy-exponential}
\end{equation}
For fixed $s$ and $0<\delta<1-1/s$, there is a constant $C_{s,\delta}>0$ such
that, for all sufficiently large $N$,
\begin{equation}
\Pr(X_N\geq t_N)
\leq C_{s,\delta}N^{-1/2}s^{-N(1-H_s(\delta))}.
\label{eq:binomial-entropy-prefactor}
\end{equation}
\end{lemma}

\begin{proof}
Using the change of variables $k=N-j$, we obtain
\begin{align*}
\Pr(X_N\geq t_N)
&=
\sum_{j=t_N}^{N}
\binom Nj
\left(\frac1s\right)^j
\left(1-\frac1s\right)^{N-j}
=
\frac1{s^N}
\sum_{k=0}^{N-t_N}\binom Nk(s-1)^k
=
\frac{V_s(N,N-t_N)}{s^N}.
\end{align*}
Since
$
N-t_N
=
N-\left(\left\lfloor(1-\delta)N\right\rfloor+1\right)
=
\lceil\delta N\rceil-1,
$
the standard estimates for the volume of a Hamming ball
\cite[Proposition~3.3.3]{GuruswamiRudraSudanECT} give
$V_s(N,N-t_N)\leq s^{NH_s(\delta)}$ and, for fixed \(s\) and
\(0<\delta\leq1-1/s\), $V_s(N,N-t_N)=s^{NH_s(\delta)+o(N)}$.
Substitution into \eqref{eq:binomial-hamming-identity} proves
\eqref{eq:binomial-entropy-upper} and
\eqref{eq:binomial-entropy-exponential}.

For \eqref{eq:binomial-entropy-prefactor}, set
$Y_N:=N-X_N\sim\operatorname{Bin}(N,1-1/s)$ and
$k_N:=N-t_N=\lceil\delta N\rceil-1$, and put $f_N:=k_N/N$ and
$p:=1-1/s$. Proposition~5.5, Eq.~(113), of
\cite{ZhuLiHayashi2022} gives
\[
\Pr(Y_N\leq k_N)
<\frac{(1-f_N)p}{\sqrt{2\pi N f_N(1-f_N)}(p-f_N)}
\exp\!\bigl(-N D(f_N\Vert p)\bigr),
\]
where $D(f\Vert p)=f\ln(f/p)+(1-f)\ln((1-f)/(1-p))$
is the binary relative entropy used in that paper. For all sufficiently
large $N$, $0<f_N\leq\delta<p$, so we have
$D(f_N\Vert p)\geq D(\delta\Vert p)$: indeed,
$\frac{d}{df}D(f\Vert p)=\ln\frac{f(1-p)}{p(1-f)}<0$ for $f<p$.
Also, since $f_N\to\delta$, the coefficient apart from $N^{-1/2}$
is bounded. To rewrite the exponent, the definition of $D$ gives
\[
D(\delta\Vert 1-1/s)
=\ln s+\delta\ln\delta-\delta\ln(s-1)
+(1-\delta)\ln(1-\delta)
=\ln s\,\bigl(1-H_s(\delta)\bigr),
\]
where the last equality follows from \eqref{eq:s-ary-entropy}. Therefore
$\exp(-N D(f_N\Vert p))\leq
s^{-N(1-H_s(\delta))}$. Together with the bounded prefactor above,
this proves \eqref{eq:binomial-entropy-prefactor}.
\end{proof}

The next observation allows us to replace the entropy associated with an
arbitrary subgroup index by the entropy at the smallest possible index. This
will be used repeatedly in the later union bounds.

\begin{lemma}[Entropy monotonicity]
\label{lem:entropy-monotonicity}
Let $q\geq2$ and $0<\delta\leq1-1/q$. Then
$H_s(\delta)\leq H_q(\delta)$ for every real $s\geq q$; for integer
$s$, this inequality is also stated in \cite[Exercise~3.6]{GuruswamiRudraSudanECT}.
\end{lemma}
\begin{proof}
Write $h(x)=-x\ln x-(1-x)\ln(1-x)$. For real $u\geq q$, set
\[
F(u):=H_u(\delta)=\frac{h(\delta)+\delta\ln(u-1)}{\ln u}.
\]
Differentiation gives
\[
F'(u)=\frac{\delta}{u(\ln u)^2}
\left(\frac{h(1-1/u)}{1-1/u}-\frac{h(\delta)}{\delta}\right).
\]
Since $(h(x)/x)'=\ln(1-x)/x^2<0$ on $(0,1)$ and
$\delta\leq1-1/q\leq1-1/u$, the term in parentheses is nonpositive.
Thus $F$ is nonincreasing for $u\geq q$, proving the claim.
\end{proof}

\section{Random Universal Codes}
\label{sec:random-univ-codes}
We begin by proving \Cref{thm:random-hadamard-not-robust-intro}, which says
that choosing a linear number of random elements from the Hadamard code does
not give rise to a good code.
\begin{theorem}
\label{thm:random-hadamard-not-robust}
    For each $k$, let $\mathbb P_k$ be the uniform probability measure on
    $\mathcal H_k$ (the universal Hadamard code as defined in \Cref{subsec:random-unic-codes-intro}). For $C>1$, put
    $
    N_k(C):=\lceil Ck\rceil,
    $
    and let $A_k(C)$ consist of $N_k(C)$ independent
    $\mathbb P_k$-distributed elements of $F_k$. Then
    \[
    \bbP\left[
    A_k(C)\text{ freely generates a free subgroup of rank }N_k(C)
    \right]
    \xrightarrow[k\to\infty]{}1.
    \]
    Consequently, in the terminology of \Cref{def:probabilistic-GV}, for
    $\cG=(F_k,\mathbb P_k)_{k\geq1}$,
    $\bbP_{\cG}(\delta,1/C)=0$ for every $0<\delta<1$. In particular,
    $\cG$ is not a $(\delta,\rho)$-GV family for any
    $0<\delta,\rho<1$.
\end{theorem}

The mechanism behind this negative result is that a typical tuple freely
generates the subgroup it spans. We isolate the combinatorial condition used
to prove this freeness.

\begin{definition}[Central tree property]
\label{def:central-tree-property}
Let $\mathbf u=(u_1,\ldots,u_N)$ be a tuple of nontrivial reduced words, and
write $\mathbf u^{\pm1}:=(u_1,u_1^{-1},\ldots,u_N,u_N^{-1})$.
Let $\operatorname{lcp}(\mathbf u^{\pm1})$ be the maximum length of a common
prefix of two distinct words in $\mathbf u^{\pm1}$, and set
$
\min(\mathbf u):=\min_{1\leq i\leq N}|u_i|.
$
We say that $\mathbf u$ has the \emph{central tree property} if
\[
2\operatorname{lcp}(\mathbf u^{\pm1})<\min(\mathbf u).
\]
\end{definition}

The relevance of the central tree property is that it prevents excessive
cancellation between consecutive words, and therefore forces freeness.

\begin{lemma}[Central tree property]
\label{lem:central-tree-free-basis}
If a tuple $\mathbf u=(u_1,\ldots,u_N)$ has the central tree property,
then it is a free basis of the subgroup
$
\langle u_1,\ldots,u_N\rangle.
$
\end{lemma}
\begin{proof}
Set $L:=\operatorname{lcp}(\mathbf u^{\pm1})$ and
$m:=\min(\mathbf u)$.
Let \(v,w\in\mathbf u^{\pm1}\) with \(vw\neq1\). The cancellation in
\(vw\) is the length of the common prefix of \(v^{-1}\) and \(w\), and
is therefore at most \(L\). Hence
\[
|vw|
\geq |v|+|w|-2L
>
\max\{|v|,|w|\},
\]
where the strict inequality follows from
\(2L<m\leq\min\{|v|,|w|\}\).

Similarly, let \(v,w,z\in\mathbf u^{\pm1}\) with \(vw\neq1\) and
\(wz\neq1\). At most \(L\) letters cancel at each of the two
junctions. Since \(|w|\geq m>2L\), these cancellations cannot remove
the entire middle word, and therefore
\[
|vwz|
\geq |v|+|w|+|z|-4L
>
|v|-|w|+|z|.
\]
Together with the assumption that the \(u_i\) are nontrivial, these
are precisely the Nielsen reduction conditions. Thus \(\mathbf u\) is
Nielsen reduced, and a Nielsen-reduced
tuple is a free basis of the subgroup it generates; see
\cite[Chapter~I, Section~2]{LyndonSchupp1977}. The conclusion also
follows directly from \cite[Proposition~1.2]{Bassino_2016}.
\end{proof}

For a random Hadamard tuple, it therefore remains to control the individual
word lengths and the overlaps between pairs of signed words. The following
estimate supplies the required overlap bound.

\begin{lemma}[Common prefixes and suffixes]
\label{lem:hadamard-overlap}
Let \(S,T\subseteq[k]\) be independent and uniformly distributed. Then,
for every \(L\geq1\) and every \(\varepsilon,\eta\in\{\pm1\}\),
\[
\bbP\!\left[
\lcp\!\left(x_S^\varepsilon,x_T^\eta\right)\geq L
\right]
\leq 3^{-L}.
\]
\end{lemma}

\begin{proof}

Let \(S,T\subseteq[k]\) be independent and uniformly distributed random
subsets, and fix $A=\{a_1<\cdots<a_r\}\subseteq[k]$.
Fix \(a_\ell\in A\). Write
\[
E_\ell
:=
\bigcap_{x\in A\setminus\{a_\ell\}}\{x\in S\cup T\}
\cap
\bigcap_{x\in[k]\setminus A}\{x\notin S\cup T\}.
\]
Then
\[
\{S\cup T=A\}
=
\{a_\ell\in S\cup T\}\cap E_\ell.
\]
The event \(E_\ell\) depends only on the membership choices at elements
\(x\neq a_\ell\), and is therefore independent of the membership choices
at \(a_\ell\). Consequently,
\[
\begin{aligned}
\bbP[a_\ell\in S\cap T\mid S\cup T=A]
&=
\frac{
\bbP[a_\ell\in S\cap T,\ E_\ell]
}{
\bbP[a_\ell\in S\cup T,\ E_\ell]
} 
=
\frac{
\bbP[a_\ell\in S\cap T]\bbP[E_\ell]
}{
\bbP[a_\ell\in S\cup T]\bbP[E_\ell]
} 
=
\frac{\bbP[a_\ell\in S\cap T]}
     {\bbP[a_\ell\in S\cup T]}.
\end{aligned}
\]
Since $\bbP[a_\ell\in S\cap T]=1/4$ and
$\bbP[a_\ell\in S\cup T]=3/4$, we obtain
$\bbP[a_\ell\in S\cap T\mid S\cup T=A]=1/3$.
If $r\geq L$, then $\lcp(x_S,x_T)\geq L$ precisely when
$a_1,\ldots,a_L\in S\cap T$, while if $r<L$ this event is impossible.
Conditioned on $S\cup T=A$, these membership choices are independent for
the different elements of $A$. Hence, for every $A\subseteq[k]$,
\[
\bbP[\lcp(x_S,x_T)\geq L\mid S\cup T=A]\leq3^{-L}.
\]
Therefore
\[
\begin{aligned}
\bbP[\lcp(x_S,x_T)\geq L]
&=\sum_A
\bbP[S\cup T=A]\,
\bbP[\lcp(x_S,x_T)\geq L\mid S\cup T=A]
\leq 3^{-L}\sum_A\bbP[S\cup T=A]
=3^{-L}.
\end{aligned}
\]
The same argument applied to the last $L$ elements of $S\cup T$ gives
$\bbP[\lcp(x_S^{-1},x_T^{-1})\geq L]\leq3^{-L}$. Moreover, for every
$B\subseteq[k]$, the words $x_B$ and $x_B^{-1}$ contain letters with
opposite signs, and hence
$\bbP[\lcp(x_S^{\pm1},x_T^{\mp1})\geq L]=0\leq3^{-L}$.
\end{proof}

\begin{proof}[Proof of \Cref{thm:random-hadamard-not-robust}]
Put $N:=N_k(C)$, write $A_k:=A_k(C)=(u_1,\ldots,u_N)$, and set
$L:=\lceil4\log_3k\rceil$.
Since \(|u_j|\sim\operatorname{Bin}(k,1/2)\) and
\(\mathbb{E}[|u_j|]=k/2\), the lower-tail Chernoff bound
\cite[Theorem~4.5(2)]{Mitzenmacher_Upfal_2005}, followed by a union bound,
gives
\[
\bbP\left[\min_{1\leq j\leq N}|u_j|<\frac{k}{3}\right]
\leq
\sum_{j=1}^N
\bbP\left[|u_j|<
\left(1-\frac13\right)\frac{k}{2}\right]
\leq
N\exp\left(-\frac{(1/3)^2(k/2)}{2}\right)
=
Ne^{-k/36}.
\]
Moreover, by \Cref{lem:hadamard-overlap},
\[
\begin{aligned}
\bbP[\lcp(A_k^{\pm1})\geq L]
&=
\bbP\left[
\exists\,1\leq i<j\leq N,\ 
\varepsilon,\eta\in\{\pm1\}:
\lcp(u_i^\varepsilon,u_j^\eta)\geq L
\right] \\
&\leq
\sum_{1\leq i<j\leq N}
\sum_{\varepsilon\in\{\pm1\}}
\bbP\left[
\lcp(u_i^\varepsilon,u_j^\varepsilon)\geq L
\right] 
\leq
2\binom N2 3^{-L}
\leq
2\binom N2 k^{-4}.
\end{aligned}
\]
Therefore, with probability at least
$1-Ne^{-k/36}-2\binom N2k^{-4}$,
every $u_j$ has length at least $k/3$, while all common prefixes of
distinct words in $\{u_1^{\pm1},\ldots,u_N^{\pm1}\}$ have length less than
$L$. Indeed, overlaps are controlled by $\lcp$, and opposite-sign words
have no common prefix.

For sufficiently large $k$, we have $k/3>2L$. Hence
$(u_1,\ldots,u_N)$ satisfies the central tree property and therefore, by
\Cref{lem:central-tree-free-basis}, freely
generates a free subgroup of rank $N$. Since $N=\lceil Ck\rceil>k$, this
subgroup cannot equal $F_k$. Thus $\langle A_k \rangle\lneqq F_k$ and
$\delta(A_k;F_k)=0$.
Finally, since $N\leq(C+1)k$,
\[
\bbP\left[\delta(A_k;F_k)=0\right]
\geq
1-(C+1)ke^{-k/36}-\frac{(C+1)^2}{k^2}\xrightarrow[k\to\infty]{}1.
\]
\end{proof}
Similarly, sampling a linear number of elements uniformly from word-metric
balls of radius at least $7$ in $F_k$ does not produce good codes.
Here and below, $B_{F_k}(R)$ denotes the radius-$R$ ball in the word metric
associated with the standard free generating set
$\{x_1,\ldots,x_k\}$.

\begin{theorem}
\label{thm:ball-not-GV}
Let $(R_k)_{k\geq1}$ be a sequence of nonnegative integers such that
$R_k\geq7$ for all sufficiently large $k$, and let $\mathbb P_k$ be the
uniform probability measure on $B_{F_k}(R_k)$. For $C>1$, put
$
N_k(C):=\lceil Ck\rceil,
$
and let $A_k(C)$ consist of $N_k(C)$ independent
$\mathbb P_k$-distributed elements of $F_k$. Then
\[
\bbP\left[
A_k(C)\text{ freely generates a free subgroup of rank }N_k(C)
\right]
\xrightarrow[k\to\infty]{}1.
\]
Consequently, in the terminology of \Cref{def:probabilistic-GV}, for
$\cG=(F_k,\mathbb P_k)_{k\geq1}$,
$\bbP_{\cG}(\delta,1/C)=0$ for every $0<\delta<1$. In particular,
$\cG$ is not a $(\delta,\rho)$-GV family for any
$0<\delta,\rho<1$.
\end{theorem}

\begin{proof}
Put $N:=N_k(C)$, write $\mathbf u:=A_k(C)=(u_1,\ldots,u_N)$, and set
$q_k:=2k-1$.
Since there are
$2kq_k^{\ell-1}$ reduced words of length $\ell\geq1$, we have
\[
|B_{F_k}(r)|
=
1+\sum_{\ell=1}^{r}2kq_k^{\ell-1}
=
\frac{kq_k^r-1}{k-1}.
\]
Since $R_k\geq7$ eventually,
\[
\bbP[|u_i|<6]
=
\frac{|B_{F_k}(5)|}{|B_{F_k}(R_k)|}
\leq
q_k^{5-R_k},
\]
and therefore
\[
\bbP\left[\min_{1\leq i\leq N}|u_i|<6\right]
\leq
Nq_k^{5-R_k}.
\]

Conditioned on $\min_i|u_i|\geq6$, two signed words arising from
distinct $u_i$ have the same length-$3$ prefix with probability
$(2kq_k^2)^{-1}\leq q_k^{-3}$. For the pair $u_i,u_i^{-1}$,
conditioned on $|u_i|=m\geq6$, fixing the first $m-3$ letters leaves
$q_k^3$ possible extensions, at most one of which satisfies
$
\lcp(u_i,u_i^{-1})\geq3.
$
Hence
\[
\bbP\left[
\lcp(\mathbf u^{\pm1})\geq3
\,\middle|\,
\min_i|u_i|\geq6
\right]
\leq
4\binom N2\frac1{2kq_k^2}+Nq_k^{-3}
\leq
\binom{2N}{2}q_k^{-3}.
\]
It follows that
\[
\begin{aligned}
\bbP\left[
2{\lcp(\mathbf u^{\pm1})}\geq\min_{1\leq i\leq N}|u_i|
\right]
&\leq
Nq_k^{5-R_k}+\binom{2N}{2}q_k^{-3} \\
&\leq
(C+1)k(2k-1)^{-2}
+
2(C+1)^2k^2(2k-1)^{-3}
=
O(k^{-1}).
\end{aligned}
\]
Thus $A_k(C)$ has the central tree property with probability tending
to $1$. By \Cref{lem:central-tree-free-basis}, it therefore freely
generates a free subgroup of rank $N_k(C)$ with probability tending
to $1$.

Since $C>1$, we have $N_k(C)>k$, so this subgroup cannot equal $F_k$.
Hence $A_k(C)$ does not generate $F_k$, and
\cite[Fact~2.1]{ChapmanDinurLubotzky} gives
$
\delta(A_k(C);F_k)=0.
$
This also proves the final assertion of the theorem.
\end{proof}

\section{Random Codes in Finite Groups}
\label{sec:finite-groups}
Throughout, let $\mathcal G=(G_i)_{i\geq1}$ be a family of nontrivial
finite groups. When no sampling measures are specified, we take
$\mathbb P_i$ to be the uniform probability measure on $G_i$. The notions
of GV and strong GV are those of \Cref{def:probabilistic-GV}; the GV rate,
range, and family parameters are defined in \Cref{def:GV-rate}. We record
the maximal-subgroup parameters used below:

\begin{definition}[Maximal-subgroup parameters]
\label{def:maximal-subgroup-parameters}
Let $\Gamma$ be a nontrivial finite group. For $s\geq2$, define
\[
m_s(\Gamma):=\#\{M<_{\max}\Gamma:[\Gamma:M]=s\},
\]
the number of maximal subgroups of $\Gamma$ of index $s$. Define also
\[
q(\Gamma):=\min\{s\geq2:m_s(\Gamma)>0\},
\qquad
\mathcal M(\Gamma):=
\max_{\substack{s\geq2\\m_s(\Gamma)>0}}
\frac{\log m_s(\Gamma)}{\log s}.
\]
Thus $q(\Gamma)$ is the smallest index of a proper subgroup, while
$\mathcal M(\Gamma)$ is its maximal-subgroup growth exponent.
\end{definition}

We introduced some of these family parameters in the introduction; we now
recall them using the precise maximal-subgroup notation above.

\begin{definition}[Parameters of a family]
\label{def:family-parameters-finite}
Recall the parameters of $\cG=(G_i)_{i\geq1}$. Put
$r_i:=\rk(G_i)$ and define
\[
q_{\cG}:=\liminf_{i\to\infty}q(G_i)
\in\{2,3,\ldots\}\cup\{\infty\},\qquad
\mu_{\cG}:=\max\left\{1,\limsup_{i\to\infty}
\frac{\mathcal M(G_i)}{r_i}\right\},\qquad
\rk_{\cG}:=\liminf_{i\to\infty}r_i\in\bbN\cup\{\infty\}.
\]
\end{definition}

\begin{definition}[Uniform polynomial maximal-subgroup growth]
\label{def:uniform-PMSG}
We say that $\cG$ has \emph{uniform polynomial maximal-subgroup growth},
or \emph{uniform PMSG}, if there is a constant $A\geq0$ such that, for
all sufficiently large $i$ and every $s\geq2$,
\[
m_s(G_i)\leq s^{A\rk(G_i)}.
\]
\end{definition}

We also require a standard parameter measuring random generation.
\begin{definition}
\label{def:prob-param}
For $k\geq1$, define
$P_\Gamma(k):=\bbP[\langle g_1,\ldots,g_k\rangle=\Gamma]$,
where $g_1,\ldots,g_k$ are independent and uniform in $\Gamma$. For
$\eta>1$, put
\[
\nu_\eta(\Gamma)
:=\min\left\{k\geq1:P_\Gamma(k)\geq\frac1\eta\right\},
\qquad
\nu(\Gamma):=\nu_e(\Gamma).
\]
Thus $\nu(\Gamma)$ is the least $k$ for which $k$ independent uniform
elements generate $\Gamma$ with probability at least $1/e$.
\end{definition}
The following lemma gives the basic estimate for the probability that
a uniformly random tuple fails to have detection probability at least
$\delta$, in terms of the maximal-subgroup growth.

\begin{lemma}
\label{lem:entropy-sensitive-failure}
Let $\Gamma$ be a nontrivial finite group, put
$r=\rk(\Gamma)$ and $q=q(\Gamma)$, and let
$0<\rho<1$ and $0<\delta<1-1/q$. Set
$C:=1/\rho$, $N(C):=\lceil Cr\rceil$, and let
$A_\Gamma(C)\in\Gamma^{N(C)}$
be a uniformly random tuple. Then
\[
\bbP[\delta(A_\Gamma(C);\Gamma)<\delta]
\leq
\sum_{s=q}^{|\Gamma|}
m_s(\Gamma)s^{-Cr(1-H_q(\delta))}.
\]
\end{lemma}

\begin{proof}
Write $A_\Gamma(C)=(X_1,\ldots,X_{N(C)})$ and put
$t:=\lceil(1-\delta)N(C)\rceil$. If $M<_{\max}\Gamma$ has index $s$,
then
\[
\bbP[X_i\in M]=\frac1s,
\qquad
|A_\Gamma(C)\cap M|
\sim\operatorname{Bin}\left(N(C),\frac1s\right).
\]
Hence, by \Cref{lem:binomial-entropy-tail},
\[
\bbP[|A_\Gamma(C)\cap M|\geq t]
\leq s^{-N(C)(1-H_s(\delta))}.
\]
If $\delta(A_\Gamma(C);\Gamma)<\delta$, then some maximal subgroup
contains at least $t$ sampled elements. Therefore,
\[
\bbP[\delta(A_\Gamma(C);\Gamma)<\delta]
\leq
\bbP\left[
\exists\,M<_{\max}\Gamma:
|A_\Gamma(C)\cap M|\geq t
\right].
\]
Therefore, a union bound gives
\[
\begin{aligned}
\bbP[\delta(A_\Gamma(C);\Gamma)<\delta]
&\leq
\sum_{s=q}^{|\Gamma|}
m_s(\Gamma)s^{-N(C)(1-H_s(\delta))}
\leq
\sum_{s=q}^{|\Gamma|}
m_s(\Gamma)s^{-Cr(1-H_q(\delta))}.
\end{aligned}
\]
Indeed, the last inequality follows from $N(C)\geq Cr$ and
$H_s(\delta)\leq H_q(\delta)$ for $s\geq q$ by
\Cref{lem:entropy-monotonicity}.
\end{proof}

We now combine this failure-probability estimate with bounds on maximal-subgroup
growth. This yields the qualitative characterization of GV families.

\begin{theorem}[GV characterization]
    \label{thm:GV-characterization}
    Let \(\cG=(G_i)_{i\geq1}\) be a family of nontrivial finite groups,
    with \(q_{\cG}\) as in \Cref{def:family-parameters}, and adopt the
    convention that \(1/q_{\cG}=0\) if \(q_{\cG}=\infty\). Then the following are equivalent:
    \begin{enumerate}
        \item The family \(\cG\) is \((\delta,\rho)\)-GV for some
        $0<\delta,\rho<1$.
    
        \item The family \(\cG\) has uniform PMSG.
    
        \item For every
        $
        0<\delta<1-1/q_{\cG},
        $
        there exists \(0<\rho<1\) such that
        \(\cG\) is \((\delta,\rho)\)-GV.
    \end{enumerate}
\end{theorem}

\begin{proof}
Suppose first that $\mathcal G$ is $(\delta,\rho)$-GV\@. Let
$N_i:=\lceil r_i/\rho\rceil$, let
$A_i=(X_{i,1},\ldots,X_{i,N_i})$ consist of independent uniform
elements of $G_i$, and put
$K_i:=\lfloor(1-\delta)N_i\rfloor+1$.

On the event $\delta(A_i;G_i)\geq\delta$, every proper subgroup of
$G_i$ contains fewer than $K_i$ elements of $A_i$. Consequently,
\[
\langle X_{i,1},\ldots,X_{i,K_i}\rangle=G_i,
\qquad
P_{G_i}(K_i)\geq
\bbP[\delta(A_i;G_i)\geq\delta].
\]
Set $\lambda:=\bbP_{\cG}(\delta,\rho)>0$
and choose $\eta>1/\lambda$. For all sufficiently large $i$,
$P_{G_i}(K_i)>1/\eta$, so $\nu_\eta(G_i)\leq K_i$. By
\cite[Lemma~2.1]{LucchiniMoscatiello},
\[
\mathcal M(G_i)
\leq\nu_\eta(G_i)+2+\log_2\eta
\leq K_i+2+\log_2\eta.
\]
Moreover,
\[
K_i\leq(1-\delta)\lceil r_i/\rho\rceil+1
\leq\frac{1-\delta}{\rho}r_i+2,
\]
and therefore
\[
\mathcal M(G_i)
\leq\frac{1-\delta}{\rho}r_i+4+\log_2\eta.
\]
Since $m_s(G_i)\leq s^{\mathcal M(G_i)}$ for every $s\geq2$,
and $r_i\geq1$, it follows that
\[
m_s(G_i)
\leq
s^{\left(\frac{1-\delta}{\rho}+4+\log_2\eta\right)r_i}.
\]
Hence the family $\mathcal G$ has uniform PMSG\@. Thus
$\textup{(1)}\Rightarrow\textup{(2)}$.

Suppose now that \textup{(2)} holds. Thus there exists $A\geq0$
such that
$
m_s(G_i)\leq s^{A\cdot r_i}
$
for all sufficiently large $i$ and every $s\geq2$. Fix
$0<\delta<1-1/q_{\cG}$. If $q_{\cG}<\infty$, put
$q=q_{\cG}$. If $q_{\cG}=\infty$, choose any sufficiently large integer
$q\geq2$ for which $\delta<1-1/q$.

Put $L:=1-H_q(\delta)>0$ and choose $0<\rho<1$ such that
$L/\rho\geq A+2$. Let $A_i$ and $N_i$ be defined as before.
By \Cref{lem:entropy-sensitive-failure},
\[
\bbP[\delta(A_i;G_i)<\delta]
\leq\sum_{s=q}^{\infty}s^{A\cdot r_i-N_iL}.
\]
Since $N_i\geq r_i/\rho$ and $r_i\geq1$,
\[
Ar_i-N_iL
\leq(A-L/\rho)r_i
\leq-2r_i
\leq-2.
\]
Consequently,
\[
\bbP[\delta(A_i;G_i)<\delta]
\leq\sum_{s=q}^{\infty}s^{-2}
\leq\sum_{s=2}^{\infty}s^{-2}
=\frac{\pi^2}{6}-1<1.
\]
Thus
\[
\bbP[\delta(A_i;G_i)\geq\delta]
\geq2-\frac{\pi^2}{6}>0
\]
for all sufficiently large $i$. Hence $\mathcal G$ is
$(\delta,\rho)$-GV, proving
$\textup{(2)}\Rightarrow\textup{(3)}$. Since the interval
$0<\delta<1-1/q_{\cG}$ is nonempty,
$\textup{(3)}\Rightarrow\textup{(1)}$ is immediate.
\end{proof}
The next two theorems refine the preceding qualitative characterization
by giving explicit bounds on the GV rate, treating separately the cases
$\rk_{\cG}<\infty$ and $\rk_{\cG}=\infty$.
\begin{theorem}[Quantitative bounds and optimal GV range for
$\rk_{\cG}=\infty$]
    \label{thm:GV-range}
    Let \(\cG=(G_i)_{i\geq1}\) be a family of nontrivial finite groups with uniform PMSG
    such that \(\rk_{\cG} = \infty\).
    We adopt the conventions \(1/q_{\cG}=0\) and
    \(H_{q_{\cG}}(\delta)=\delta\) if \(q_{\cG}=\infty\). Then:
    \begin{enumerate}
        \item\label{thm:GV-range-part1-finite-groups}
        For every \(0<\delta<1-1/q_{\cG}\),
        \[
        \frac{1-H_{q_{\cG}}(\delta)}{\mu_{\cG}}
        \leq R_{\cG}(\delta)
        \leq
        \frac{1-\delta}{\mu_{\cG}}.
        \]
        Moreover, if
        $
        0<\rho<
        \frac{1-H_{q_{\cG}}(\delta)}{\mu_{\cG}}
        $
        then
        $
        \bbP_{\cG}(\delta,\rho)=1.
        $
    
        \item If \(q_{\cG}=\infty\), then, for every \(0<\delta<1\),
        $
        R_{\cG}(\delta)
        =
        \frac{1-\delta}{\mu_{\cG}}.
        $
        \item If \(q_{\cG}<\infty\), then, for every
        \(1-1/q_{\cG}<\delta<1\),
        $
        R_{\cG}(\delta)=0.
        $
    \end{enumerate}
\end{theorem}
\begin{theorem}[Quantitative GV bounds for $\rk_{\cG}<\infty$]
\label{thm:GV-range-finite-d}
Let $\cG=(G_i)_{i\geq1}$ be a family of nontrivial finite groups with
uniform PMSG and $\rk_{\cG}<\infty$.
Then:
\begin{enumerate}
    \item If $q_{\cG}<\infty$, then, for every
    $0<\delta<1-1/q_{\cG}$, we have
    $
    R_{\cG}(\delta)
    \geq
    \frac{1-H_{q_{\cG}}(\delta)}
    {\mu_{\cG}+2/\rk_{\cG}}.
    $
    \item If $q_{\cG}=\infty$, then, for every $0<\delta<1$, we have
    $
    R_{\cG}(\delta)
    \geq
    \frac{1-\delta}{\mu_{\cG}+1/\rk_{\cG}}.
    $

    Moreover, if
    $
    0<\rho<\frac{1-\delta}{\mu_{\cG}+1/\rk_{\cG}}
    $
    then
    $
    \bbP_{\cG}(\delta,\rho)=1.
    $
\end{enumerate}
\end{theorem}
The proofs of the two preceding theorems rely on the same probability
estimate, so we prove them together.

\begin{proof}[Proof of
\Cref{thm:GV-range-finite-d,thm:GV-range}]
Put $q_i:=q(G_i)$ and $q_\cG=\liminf_{i\to\infty}q_i$. Since $\cG$ has uniform PMSG,
$\mu_{\cG}<\infty$, and, for every $\varepsilon>0$,
$
\mathcal M(G_i)\leq(\mu_{\cG}+\varepsilon)r_i
$
for all sufficiently large $i$.

Throughout the proof, whenever $0<\rho<1$ is fixed, we write
$N_i:=\lceil r_i/\rho\rceil$ and let $A_i$ consist of $N_i$
independent uniform elements of $G_i$. For
$\delta<1-1/q_i$, put $L_i:=1-H_{q_i}(\delta)$. Since
$
m_s(G_i)\leq s^{\mathcal M(G_i)}
$
for every $s\geq2$, \Cref{lem:entropy-sensitive-failure} gives
\begin{equation}
\bbP_i[\delta(A_i;G_i)<\delta]
\leq
\sum_{s=q_i}^{|G_i|}m_s(G_i)s^{-\frac{r_i}{\rho}L_i}
\leq
\sum_{s=q_i}^{\infty}
s^{\mathcal M(G_i)-\frac{r_i}{\rho}L_i}.
\label{failure_bound}
\end{equation}

For the lower-bound arguments below, fix
$0<\delta<1-1/q_{\cG}$ and put
$L:=1-H_{q_{\cG}}(\delta)>0$, with
$H_\infty(\delta)=\delta$.

\medskip
\noindent
\emph{Proof of \Cref{thm:GV-range-finite-d}.}
Suppose that $\rk_{\cG}<\infty$. Since the $r_i$ are integers,
$r_i\geq \rk_{\cG}$ for all sufficiently large $i$.

For part \textup{(1)}, suppose that $q_{\cG}<\infty$ and choose
$0<\rho<\frac{L}{\mu_{\cG}+2/\rk_{\cG}}$. Choose
$\varepsilon>0$ such that
$L/\rho\geq\mu_{\cG}+2/\rk_{\cG}+2\varepsilon$.
Since $q_i\geq q_{\cG}$ eventually, we have $L_i\geq L$ by
\Cref{lem:entropy-monotonicity}. Therefore, for all sufficiently large
$i$,
\[
\mathcal M(G_i)-\frac{r_i}{\rho}L_i
\leq
(\mu_{\cG}+\varepsilon)r_i-\frac{r_i}{\rho}L
=
\left(\mu_{\cG}+\varepsilon-\frac{L}{\rho}\right)r_i
\leq
-\left(\frac{2}{\rk_{\cG}}+\varepsilon\right)r_i
\leq-2.
\]
Hence, by \eqref{failure_bound},
\[
\bbP_i[\delta(A_i;G_i)<\delta]
\leq
\sum_{s=q_i}^{\infty}s^{-2}
\leq
\sum_{s=2}^{\infty}s^{-2}
=
\frac{\pi^2}{6}-1.
\]
Consequently, $\bbP_{\cG}(\delta,\rho)\geq2-\pi^2/6>0$.
Thus $\cG$ is $(\delta,\rho)$-GV for every
$
0<\rho<L/(\mu_{\cG}+2/\rk_{\cG}),
$
and hence
\[
R_{\cG}(\delta)
\geq
\frac{1-H_{q_{\cG}}(\delta)}
{\mu_{\cG}+2/\rk_{\cG}}.
\]

For part \textup{(2)}, suppose that $q_{\cG}=\infty$. Then
$L=1-\delta$. Choose
$0<\rho<\frac{L}{\mu_{\cG}+1/\rk_{\cG}}$.
Since $H_{q_i}(\delta)\to\delta$, we have
$
L_i/\rho\to L/\rho
>
\mu_{\cG}+1/\rk_{\cG}.
$
Thus there exists $\varepsilon>0$ such that, for all sufficiently
large $i$,
\[
\frac{L_i}{\rho}
\geq
\mu_{\cG}+\frac1{\rk_{\cG}}+2\varepsilon
\qquad\text{and}\qquad
\mathcal M(G_i)\leq(\mu_{\cG}+\varepsilon)r_i.
\]
It follows that
\[
\mathcal M(G_i)-\frac{r_i}{\rho}L_i
\leq
\left(\mu_{\cG}+\varepsilon-\frac{L_i}{\rho}\right)r_i
\leq
-\left(\frac1{\rk_{\cG}}+\varepsilon\right)r_i
\leq
-1-\varepsilon \rk_{\cG}.
\]
Hence, by \eqref{failure_bound},
\[
\bbP_i[\delta(A_i;G_i)<\delta]
\leq
\sum_{s=q_i}^{\infty}s^{-1-\varepsilon \rk_{\cG}}
\leq
\int_{q_i-1}^{\infty}x^{-1-\varepsilon \rk_{\cG}}\,dx
=
\frac{(q_i-1)^{-\varepsilon \rk_{\cG}}}
{\varepsilon \rk_{\cG}}
\xrightarrow[i\to\infty]{}0.
\]
Thus $\bbP_{\cG}(\delta,\rho)=1$, and hence
\[
R_{\cG}(\delta)
\geq
\frac{1-\delta}{\mu_{\cG}+1/\rk_{\cG}}.
\]

\medskip
\noindent
\emph{Proof of \Cref{thm:GV-range}.}
Suppose that $\rk_{\cG}=\infty$. For part \textup{(1)}, let
$\delta$ and $L$ be as above, and choose
$0<\rho<L/\mu_{\cG}$. Choose $\varepsilon>0$ such that
$L/\rho\geq\mu_{\cG}+3\varepsilon$.
If $q_{\cG}<\infty$, then $q_i\geq q_{\cG}$ eventually, so
$L_i\geq L$ by \Cref{lem:entropy-monotonicity}. If
$q_{\cG}=\infty$, then $q_i\to\infty$ and $L_i\to L$.
Thus, in either case, $L_i/\rho\geq\mu_{\cG}+2\varepsilon$ for all
sufficiently large $i$.
Together with
$\mathcal M(G_i)\leq(\mu_{\cG}+\varepsilon)r_i$, this gives
$\mathcal M(G_i)-r_iL_i/\rho\leq-\varepsilon r_i$.
Since $\varepsilon r_i\geq2$ eventually, \eqref{failure_bound} gives
\[
\bbP_i[\delta(A_i;G_i)<\delta]
\leq
\sum_{s=q_i}^{\infty}s^{-\varepsilon r_i}
\leq
2^{-(\varepsilon r_i-2)}
\sum_{s=2}^{\infty}s^{-2}
\xrightarrow[i\to\infty]{}0.
\]
This proves both the lower bound and the convergence assertion in
part \textup{(1)}.

We next prove the upper bound in part \textup{(1)}. Suppose that
$\cG$ is $(\delta,\rho)$-GV\@. The proof of
\Cref{thm:GV-characterization} gives
\[
\mathcal M(G_i)
\leq
\frac{1-\delta}{\rho}r_i+O(1).
\]
Moreover, putting
$
K_i:=\lfloor(1-\delta)\lceil r_i/\rho\rceil\rfloor+1,
$
the same proof gives
$P_{G_i}(K_i)\geq\bbP_i[\delta(A_i;G_i)\geq\delta]$.
The right-hand side is positive for all sufficiently large $i$, and
hence $P_{G_i}(K_i)>0$. Thus $G_i$ admits a generating $K_i$-tuple, so
$r_i\leq K_i$. Consequently,
\[
r_i
\leq
\left\lfloor(1-\delta)
\left\lceil\frac{r_i}{\rho}\right\rceil\right\rfloor+1
\leq
\frac{1-\delta}{\rho}r_i+2.
\]
Since $r_i\to\infty$, dividing by $r_i$ and passing to the limit yields
$
1\leq(1-\delta)/\rho.
$
Dividing the preceding bound by $r_i$ and
taking the limsup yields
$
\mu_{\cG}\leq(1-\delta)/\rho,
$
and hence
$
\rho\leq(1-\delta)/\mu_{\cG}.
$
Therefore
\[
R_{\cG}(\delta)
\leq
\frac{1-\delta}{\mu_{\cG}}.
\]

For part \textup{(2)}, suppose that $q_{\cG}=\infty$. Then
$H_{q_{\cG}}(\delta)=\delta$, so the lower and upper bounds in
part \textup{(1)} coincide. Hence, for every $0<\delta<1$,
\[
R_{\cG}(\delta)
=
\frac{1-\delta}{\mu_{\cG}}.
\]

Finally, we prove part \textup{(3)}. Suppose that
$q_{\cG}<\infty$. There is a subsequence along which
$q_i=q_{\cG}$. For each $i$ in this subsequence, choose a maximal
subgroup $M_i<G_i$ of index $q_{\cG}$.

Fix $1-1/q_{\cG}<\delta<1$ and $0<\rho<1$, and let $N_i$ and
$A_i$ be as above. Put
$
X_i:=|A_i\setminus M_i|
$
and
$
p:=1-1/q_{\cG}.
$
Then
\[
X_i\sim\operatorname{Bin}(N_i,p),
\qquad
\delta(A_i;G_i)\leq\frac{X_i}{N_i},
\]
and
\[
\mathbb{E}\left[\frac{X_i}{N_i}\right]=p,
\qquad
\operatorname{Var}\left(\frac{X_i}{N_i}\right)
=
\frac{p(1-p)}{N_i}.
\]
Put $\varepsilon:=\delta-p>0$. Since
$
\delta(A_i;G_i)\leq X_i/N_i,
$
we have
$
\{\delta(A_i;G_i)\geq\delta\}
\subseteq
\left\{\frac{X_i}{N_i}\geq p+\varepsilon\right\}.
$
Therefore, by Chebyshev's inequality,
\[
\bbP_i[\delta(A_i;G_i)\geq\delta]
\leq
\bbP_i\left[\frac{X_i}{N_i}-p\geq\varepsilon\right]
\leq
\bbP_i\left[\left|\frac{X_i}{N_i}-p\right|
\geq\varepsilon\right]
\leq
\frac{p(1-p)}{N_i\varepsilon^2}
\xrightarrow[i\to\infty]{}0,
\]
since $N_i\geq r_i/\rho\to\infty$. Thus the success probability tends
to zero along a subsequence, so
$\bbP_{\cG}(\delta,\rho)=0$ and $\cG$ is not $(\delta,\rho)$-GV for
any $0<\rho<1$. Therefore
$
R_{\cG}(\delta)=0.
$
\end{proof}

The GV characterization in \Cref{thm:GV-characterization} determines when
some nontrivial GV parameters exist. We next determine when the entire
distance range $0<\delta<1$ is attainable, that is, when $\cG$ is strong GV
in the sense of \Cref{def:strong-GV}.

\begin{proposition}[Characterization of strong GV]
\label{prop:strong-GV}
Let \(\cG=(G_i)_{i\in\bbN}\) be a family of nontrivial finite groups, and
let $q_\cG$ be as in \Cref{def:family-parameters-finite}. Then:
\begin{enumerate}
    \item If \(q_{\cG}=\infty\) and \(\cG\) has uniform PMSG, then
    \(\cG\) is strong GV.

    \item If $\rk_{\cG}=\infty$, then
    \[
    \cG\text{ is strong GV}
    \quad\Longleftrightarrow\quad
    \begin{cases}
        q_{\cG}=\infty,\\
        \cG\text{ has uniform PMSG}.
    \end{cases}
    \]
\end{enumerate}
\end{proposition}

\begin{proof}
The first assertion follows immediately from
\Cref{thm:GV-characterization}, since \(q_{\cG}=\infty\) means that the
range appearing there is \(0<\delta<1\).

For the second assertion, the reverse implication follows from the first.
Conversely, suppose that \(\cG\) is strong GV\@. In particular, it is
\((\delta,\rho)\)-GV for some \(0<\delta<1\) and \(0<\rho<1\), so
\Cref{thm:GV-characterization} implies that \(\cG\) has uniform PMSG.

It remains to show that \(q_{\cG}=\infty\). Otherwise,
\(q_{\cG}<\infty\), and we may choose
$
1-1/{q_{\cG}}<\delta<1.
$
Since \(\rk(G_i)\to\infty\), \Cref{thm:GV-range} gives
\(R_{\cG}(\delta)=0\), contradicting the assumption that \(\cG\) is
strong GV\@. Therefore \(q_{\cG}=\infty\).
\end{proof}

The assumption $\rk_{\cG}=\infty$ in the converse direction is essential. The
following examples illustrate the two different behaviors that can occur when
$\rk_{\cG}<\infty$.

\begin{remark}
\label{rem:finite-d-strong}
Uniform PMSG together with $\rk_{\cG}<\infty$ and $q_{\cG}<\infty$ does
not determine whether $\cG$ is strong GV\@. Indeed, fix a prime $p$ and
an integer $r\geq2$.

For a strong-GV example, consider
$
\cG_1=(C_{p^i})_{i\geq1}.
$
Then $\rk_{\cG_1}=1$, $q_{\cG_1}=p$, and $\cG_1$ has uniform PMSG\@.
Given $0<\delta<1$, put
$
N_\delta:=\lceil(1-\delta)^{-1}\rceil
$
and
$
\rho_\delta:=1/N_\delta\in(0,1),
$
and let $A_i$ consist of $N_\delta$ independent uniform elements of
$C_{p^i}$. If $M_i$ is the unique maximal subgroup of $C_{p^i}$, then
$\bbP_{\cG_1}(\delta,\rho_\delta)
\geq(1-1/p)^{N_\delta}>0$.
Thus $\cG_1$ is strong GV.

For a non-strong-GV example, let $\cG_2=(G_i)_{i\geq1}$ be any family
of finite $p$-groups with $\rk(G_i)=r$. By the Burnside basis theorem,
$G_i/\Phi(G_i)\cong C_p^r$.
Since $\Phi(G_i)$ is contained in every maximal subgroup of $G_i$, the
correspondence theorem identifies the maximal subgroups of $G_i$ with
the hyperplanes of $C_p^r$. Hence every maximal subgroup has index $p$,
and their number is $m_p(G_i)=(p^r-1)/(p-1)$.
Thus $\rk_{\cG_2}=r$, $q_{\cG_2}=p$, and $\cG_2$ has uniform PMSG.

Fix $0<\rho<1$, put $N:=\lceil r/\rho\rceil$, and let $A_i$ consist
of $N$ independent uniform elements of $G_i$. Denote its image in
$G_i/\Phi(G_i)$ by $\overline A_i$. By
\cite[Fact~2.1]{ChapmanDinurLubotzky},
$
\delta(A_i;G_i)=\delta(\overline A_i;C_p^r).
$
If $\overline A_i$ does not generate $C_p^r$, this quantity is zero.
Otherwise, by \cite[\S2.2]{ChapmanDinurLubotzky},
$\overline A_i$ determines a linear $[N,r,\Delta_i]_p$ code whose
relative minimum distance is $\Delta_i/N$. The Plotkin bound
\cite[Theorem~4.4.1]{GuruswamiRudraSudanECT} gives
\[
\delta(A_i;G_i)
=
\frac{\Delta_i}{N}
\leq
\frac{p^{r-1}(p-1)}{p^r-1}<1.
\]
Consequently, for every
$
p^{r-1}(p-1)/(p^r-1)<\delta<1
$
and every $0<\rho<1$,
$
\bbP_{\cG_2}(\delta,\rho)=0.
$
Thus $\cG_2$ is not strong GV, showing that the assumption
$\rk_{\cG}=\infty$ in part \textup{(2)} of
\Cref{prop:strong-GV} is essential.
\end{remark}
\section{Examples}
\label{sec:examples}
We now turn to the proof of \Cref{thm:simple-power-probabilistic-profile-intro}, stated in the introduction,
which determines the exact probabilistic GV rate for direct powers of a
fixed finite simple group. We begin by counting the maximal
subgroups of such direct powers.
\begin{lemma}
\label{lem:maximal-subgroups-simple-power}
Let \(H\) be a nonabelian finite simple group and let \(m\geq1\).
Every maximal subgroup of \(H^m\) is of exactly one of the following
types:
\begin{enumerate}[(i)]
    \item\label{item:product-type}
    a product-type subgroup
    $
    H^{j-1}\times M\times H^{m-j},
    $
    where \(j\in[m]\) and \(M\) is a maximal subgroup of \(H\);

    \item\label{item:diagonal-type}
    a diagonal-type subgroup
    $
    D_{j,k,\alpha}
    =
    \left\{
    (h_1,\ldots,h_m)\in H^m:
    h_k=\alpha(h_j)
    \right\},
    $
    where \(1\leq j<k\leq m\) and \(\alpha\in\Aut(H)\).
\end{enumerate}
Consequently, for every \(s\geq2\),
\[
m_s(H^m)
=
m\cdot m_s(H)
+
\mathbf 1_{\{s=|H|\}}
\binom{m}{2}|\Aut(H)|.
\]
\end{lemma}

\begin{proof}
By \cite[Lemma~1.3]{Thevenaz1997}, the maximal subgroups of
\(H\times H\) are precisely
$M\times H,\quad H\times M,$
where \(M<_{\max}H\), and
$
\{(h,\alpha(h)):h\in H\},
$
where \(\alpha\in\Aut(H)\).

Since \(H\) is nonabelian simple, it is perfect, and hence so is
\(H^m\). In particular, \(H^m\) has no normal subgroup of prime index.
	Therefore, \cite[Proposition~1.6]{Thevenaz1997} shows that every
	maximal subgroup of \(H^m\) is the inverse image of one of the above
	subgroups under a coordinate projection. These give precisely
	the product-type and diagonal-type subgroups in the statement.

There are \(m\cdot m_s(H)\) product-type maximal subgroups of index \(s\).
The diagonal-type subgroups all have index \(|H|\), and there is one
for each pair \(1\leq j<k\leq m\) and each
\(\alpha\in\Aut(H)\). Hence
\[
m_s(H^m)
=
m\cdot m_s(H)
+
\mathbf 1_{\{s=|H|\}}
\binom{m}{2}|\Aut(H)|.
\]
\end{proof}
The next theorem extends the sharp GV threshold for random
$p$-ary linear codes, corresponding to $H=C_p$, to direct powers of
an arbitrary finite simple group.
\begin{theorem}
\label{thm:simple-power-probabilistic-profile}
Let $H$ be a fixed finite simple group and consider the family
$\cG=\cG_H:=(H^m)_{m\in\bbN}$. For every
$0<\delta<1$,
\[
R_{\cG}(\delta)
    =
    \begin{cases}
    \log_{|H|}q(H)\bigl(1-H_{q(H)}(\delta)\bigr),
        & 0<\delta<1-\dfrac{1}{q(H)},\\[4pt]
    0,
        & 1-\dfrac{1}{q(H)}\leq\delta<1.
    \end{cases}
    \qquad
    \Range(\cG)=1-\frac{1}{q(H)}.
\]
For every $0<\delta<1-1/q(H)$ and every $0<\rho<1$,
\[
\bbP_{\cG}(\delta,\rho)
=
\begin{cases}
1, & \rho\leq R_{\cG}(\delta),\\[2mm]
0, & \rho>R_{\cG}(\delta).
\end{cases}
\]
\end{theorem}

\begin{proof}
Fix $0<\rho<1$, put $C:=1/\rho$, and set
$N_m(C):=\lceil C\,\rk(H^m)\rceil$.
Let $A_m(C)$ be a tuple of $N_m(C)$ independent uniformly distributed
elements of $H^m$.
If $H$ is abelian, then $H\cong C_p$ for some prime $p$. Identifying
$H^m$ with $\bbF_p^m$ and writing the elements of $A_m(C)$ as the rows
of a random matrix $X\in\bbF_p^{N_m(C)\times m}$, we have
\[
\delta\bigl(A_m(C);H^m\bigr)
=
\frac{1}{N_m(C)}
\min_{a\neq0}\operatorname{wt}(Xa),
\]
since the maximal subgroups of $\bbF_p^m$ are its hyperplanes. This is
the normalized minimum distance of the associated random $p$-ary linear
code. For $0<\delta<1-1/p$, by
\cite[Corollary~3.2, Remark~3.3, and Lemma~2.1]{Fan2012CumulativeDE},
$R_{\cG}(\delta)=1-H_p(\delta)$, and
\[
\bbP_{\cG}(\delta,\rho)
=
\begin{cases}
1, & \rho\leq R_{\cG}(\delta),\\[2mm]
0, & \rho>R_{\cG}(\delta).
\end{cases}
\]
Since $q(H)=|H|=p$, in this range we obtain
\[
R_{\cG}(\delta)
=
\log_{|H|}q(H)
\bigl(1-H_{q(H)}(\delta)\bigr).
\]
In this case, $q_{\cG}=p$, $\rk_{\cG}=\infty$, and $\cG$ has uniform
PMSG. The quantitative characterization gives $R_{\cG}(\delta)=0$ for
$\delta>1-1/p$ and $\Range(\cG)=1-1/p$.
At $\delta_*=1-1/p$, the corresponding random matrix defines a
\(p\)-ary linear code of length $N_m(C)$ and dimension $m$ whenever
$\delta(A_m(C);H^m)\geq\delta_*$. If this event occurred for infinitely
many $m$, these codes would have relative distance at least $\delta_*$
and rates $m/N_m(C)\to\rho>0$, contradicting \Cref{thm:plotkin}. Hence
$\bbP_{\cG}(\delta_*,\rho)=0$ for every $0<\rho<1$, and
$R_{\cG}(\delta_*)=0$.
We may therefore assume that $H$ is nonabelian.
Put $r_m:=\rk(H^m)$ and
$t_m:=\lfloor(1-\delta)N_m(C)\rfloor+1$. By
\cite[Proposition~3.4]{Thevenaz1997},
$
r_m=\log_{|H|}m+O(1).
$
Therefore,
$\log_{|H|}m/N_m(C)\to1/C=\rho$.
By \Cref{lem:maximal-subgroups-simple-power},
\[
m_s(H^m)
=
m\cdot m_s(H)
+
\mathbf 1_{\{s=|H|\}}\binom m2|\Aut(H)|.
\]
It follows that
\begin{align*}
\mathcal M(H^m)
&=
\max_{\substack{s\geq2\\m_s(H^m)>0}}
\frac{\log\left(m\,m_s(H)+
\mathbf 1_{\{s=|H|\}}\binom m2|\Aut(H)|\right)}{\log s}\\
&=
\max\left\{
\max_{\substack{s\geq2\\m_s(H)>0}}
\frac{\log\bigl(m\,m_s(H)\bigr)}{\log s},
\frac{\log\left(\binom m2|\Aut(H)|\right)}{\log|H|}
\right\}+O(1)\\
&=
\max\left\{\log_{q(H)}m,2\log_{|H|}m\right\}
+O(1).
\end{align*}
Consequently,
\[
\mu_{\cG}
=
\limsup_{m\to\infty}
\frac{\mathcal M(H^m)}{\rk(H^m)}
=
\max\left\{\log_{q(H)}|H|,2\right\}.
\]
We also have $q(H^m)=q(H)$. Indeed, placing a subgroup of
index $q(H)$ in one coordinate gives
$q(H^m)\leq q(H)$. Conversely, let $K<H^m$. Since the direct
factors generate $H^m$, some direct factor $H_j$ is not contained in
$K$. Hence $H_j\cap K<H_j$, and
$[H^m:K]\geq[H_j:H_j\cap K]\geq q(H)$.
Therefore $q(H^m)\geq q(H)$, proving the equality. For
$1-1/q(H)<\delta<1$, \Cref{thm:GV-range} gives
$R_{\cG}(\delta)=0$. We may therefore assume
$0<\delta\leq1-1/q(H)$. In this range, for
$0<\delta<1-1/q(H)$, \Cref{thm:GV-range} gives
\[
R_{\cG}(\delta)
\geq
\frac{1-H_{q(H)}(\delta)}{\mu_{\cG}}
=
\min\left\{
\log_{|H|}q(H)\bigl(1-H_{q(H)}(\delta)\bigr),
\frac{1-H_{q(H)}(\delta)}2
\right\}.
\]
By \cite[Theorem~A(4)]{BallesterBolinchesEtAl2022},
$q(H)^2<|H|$, and hence
$\log_{|H|}q(H)<1/2$.  Therefore,
\[
\log_{|H|}q(H)\bigl(1-H_{q(H)}(\delta)\bigr)
<
\frac{1-H_{q(H)}(\delta)}2.
\]
Hence, for $0<\delta<1-1/q(H)$,
\[
R_{\cG}(\delta)
\geq
\log_{|H|}q(H)\bigl(1-H_{q(H)}(\delta)\bigr).
\]
At $\delta=1-1/q(H)$, we have $1-H_{q(H)}(\delta)=0$, and the
trivial bound $R_{\cG}(\delta)\geq0$ supplies the lower bound. We now
prove the reverse inequality for
$0<\delta\leq1-1/q(H)$. Represent the elements of $A_m(C)$
as the rows of a random matrix
\[
X=(x_{ij})\in H^{N_m(C)\times m}.
\]
Suppose that
\[
\rho>
\log_{|H|}q(H)
\bigl(1-H_{q(H)}(\delta)\bigr).
\]
Fix $M<_{\max}H$ of index $q(H)$ and, for each $j\in[m]$, put
\[
Y_{m,j}
:=
\#\{i\in[N_m(C)]:x_{ij}\in M\},
\qquad
E_j:=\{Y_{m,j}\geq t_m\}.
\]
Since
$
\bbP[x_{ij}\in M]=1/q(H),
$
we have
$
Y_{m,j}
\sim
\operatorname{Bin}\left(N_m(C),\frac1{q(H)}\right).
$
Moreover, the columns of $X$ are independent, so the events
$E_1,\ldots,E_m$ are independent. Hence, by
\Cref{eq:binomial-entropy-exponential},
\[
p_m:=\bbP(E_j)
=
q(H)^{-N_m(C)(1-H_{q(H)}(\delta))
+o(N_m(C))}.
\]
Since $\log_{|H|}m/N_m(C)\xrightarrow[m\to\infty]{}\rho$, we have
\[
\frac1{N_m(C)}\log_{|H|}(mp_m)
\xrightarrow[m\to\infty]{}
\rho-
\log_{|H|}q(H)\bigl(1-H_{q(H)}(\delta)\bigr)
>0.
\]
Thus $mp_m\xrightarrow[m\to\infty]{}\infty$, and hence
\[
\bbP\left[\bigcup_{j=1}^{m}E_j\right]
=
1-(1-p_m)^m
\geq
1-e^{-mp_m}
\xrightarrow[m\to\infty]{}1.
\]

Whenever some $E_j$ occurs, at least $t_m>(1-\delta)N_m(C)$ rows of $X$
belong to the maximal subgroup
$
H^{j-1}\times M\times H^{m-j}<H^m.
$
Consequently, $\bbP_{\cG}(\delta,\rho)=0$. It follows that
\[
R_{\cG}(\delta)
\leq
\log_{|H|}q(H)\bigl(1-H_{q(H)}(\delta)\bigr).
\]
For $0<\delta<1-1/q(H)$, combining the two inequalities proves the
rate formula. At $\delta=1-1/q(H)$, the right-hand side is zero, so this
same reverse inequality proves $R_{\cG}(\delta)=0$. We may therefore assume
$0<\delta<1-1/q(H)$ in the following threshold case, so that
$\rho_0>0$.

We now treat the threshold value
$\rho_0:=\log_{|H|}q(H)(1-H_{q(H)}(\delta))$.
For brevity, put $N_m:=N_m(1/\rho_0)$. Since
$\rk(H^m)=\log_{|H|}m+O(1)$ and
$N_m=\lceil\rk(H^m)/\rho_0\rceil$, we have
\[
\log_{|H|}m=N_m\rho_0+O(1).
\]
Since every maximal subgroup of $H^m$ has index at most $|H|$, the
possible indices $s$ lie in the fixed finite set
$\{q(H),\ldots,|H|\}$.
If $\delta(A_m(1/\rho_0);H^m)<\delta$, then some maximal subgroup
contains at least $\lfloor(1-\delta)N_m\rfloor+1$ rows. For each maximal subgroup of index $s$,
the binomial estimate \Cref{eq:binomial-entropy-prefactor} gives a
failure probability at most
$K_{|H|,\delta}N_m^{-1/2}s^{-N_m(1-H_s(\delta))}$, where the constant
is uniform over this finite set of indices. A union bound, followed by
substitution of the formula for $m_s(H^m)$, gives
\[
\begin{aligned}
&\bbP\bigl[\delta(A_m(1/\rho_0);H^m)<\delta\bigr]\\
&\quad\leq K_{|H|,\delta}N_m^{-1/2}
\sum_{s=q(H)}^{|H|}
\left(m\,m_s(H)+\mathbf 1_{\{s=|H|\}}
\binom m2|\Aut(H)|\right)s^{-N_m(1-H_s(\delta))}\\
&\quad=K_{|H|,\delta}N_m^{-1/2}
\underbrace{\left(
\begin{aligned}
&\sum_{s=q(H)}^{|H|}m\,m_s(H)s^{-N_m(1-H_s(\delta))}
&+\binom m2|\Aut(H)|\,|H|^{-N_m(1-H_{|H|}(\delta))}
\end{aligned}
\right)}_{=:B_m}.
\end{aligned}
\]
For $s=q(H)$, the corresponding summand is bounded because
\[
m q(H)^{-N_m(1-H_{q(H)}(\delta))}
=|H|^{\log_{|H|}m-N_m\rho_0}=O(1).
\]
Put
\[
\phi(s):=\log_{|H|}s\,(1-H_s(\delta)).
\]
By \Cref{lem:entropy-monotonicity} and positivity of
$1-H_s(\delta)$, $\phi$ is strictly increasing in $s$. In particular, for each
$s>q(H)$ we have $\phi(s)>\phi(q(H))=\rho_0$. Hence
\[
\begin{aligned}
m\,s^{-N_m(1-H_s(\delta))}
&=|H|^{\log_{|H|}m-N_m\phi(s)}
&=|H|^{N_m(\rho_0-\phi(s))+O(1)}
\longrightarrow0.
\end{aligned}
\]
The convergence is exponential since $\phi(s)>\rho_0$. Finally,
\[
\binom m2|H|^{-N_m(1-H_{|H|}(\delta))}
=O\!\left(
|H|^{N_m(2\rho_0-(1-H_{|H|}(\delta)))}
\right)\longrightarrow0,
\]
because $q(H)^2<|H|$ and, by the strict increase of $\phi$,
\[
1-H_{|H|}(\delta)
=\phi(|H|)
>\phi(q(H)^2)
\geq\log_{|H|}(q(H)^2)(1-H_{q(H)}(\delta))
=2\rho_0.
\]
The estimates above show that $B_m=O(1)$. Hence the factor
$N_m^{-1/2}$ forces the failure probability to zero. Therefore
$\bbP_{\cG}(\delta,\rho_0)=1$.
For nonabelian $H$, the preceding computation gives $q_{\cG}=q(H)$ and
$\mu_{\cG}<\infty$, while $\rk_{\cG}=\infty$. Thus
\Cref{thm:GV-range} gives $R_{\cG}(\delta)=0$ for
$1-1/q(H)<\delta<1$. The exact rate formula is positive for
$0<\delta<1-1/q(H)$.
Therefore
\[
\Range(\cG)=1-\frac{1}{q(H)}.
\]
\end{proof}

The preceding computation is exact for powers of a fixed simple group. The
general rate bound also gives useful conclusions for families whose orders
grow moderately relative to their generating ranks.

\begin{corollary}
\label{cor:probabilistic-qary-gv}
Let $\cG=(G_i)_{i\geq1}$ be a family of nontrivial finite groups. Put $r_i=\rk(G_i)$ and assume that $\rk_{\cG}=\infty$ and $\log\log|G_i|=o(r_i)$.
Suppose that $q_{\cG}<\infty$. Then, for every
$0<\delta<1-1/q_{\cG}$,
\[
R_{\cG}(\delta)
\geq
1-H_{q_{\cG}}(\delta).
\]
Moreover,
$
\bbP_{\cG}(\delta,\rho)=1
$
for every
$0<\rho<1-H_{q_{\cG}}(\delta)$.
\end{corollary}

\begin{proof}
By \cite[Theorem~2.1]{LubotzkyRandomGeneration},
$\mathcal M(G_i)\leq r_i+2\log\log|G_i|+2$.
Hence $\mu_{\cG}=1$.
Therefore, \Cref{thm:GV-range} gives
\[
R_{\cG}(\delta)
\geq
\frac{1-H_{q_{\cG}}(\delta)}{\mu_{\cG}}
=
1-H_{q_{\cG}}(\delta).
\]
The assertion about $\bbP_{\cG}$ follows from the convergence statement in
\Cref{thm:GV-range}.
\end{proof}

Direct powers of a fixed nonperfect group satisfy these growth hypotheses
automatically, giving the following concrete application.

\begin{corollary}[Direct powers of nonperfect groups]
\label{cor:direct-powers-nonperfect}
Let $H$ be a fixed finite nonperfect group, and let
$\cG:=\cG_H=(H^m)_{m\geq1}$. Then $q_{\cG}=q(H)$ and, for every
$0<\delta<1-1/q_{\cG}$,
\[
R_{\cG}(\delta)
\geq
1-H_{q_{\cG}}(\delta).
\]
Moreover,
$
\bbP_{\cG}(\delta,\rho)=1
$
for every
$0<\rho<1-H_{q_{\cG}}(\delta)$.
\end{corollary}

\begin{proof}
Since $H$ is nonperfect, it admits a quotient $H\twoheadrightarrow C_p$
for some prime $p$. Hence $H^m\twoheadrightarrow C_p^m$, and therefore
$m\leq \rk(H^m)\leq m\,\rk(H)$. Thus $\rk(H^m)\to\infty$ and
$\log_2\log_2|H^m|=\log_2m+O(1)=o(\rk(H^m))$.
As in the proof of \Cref{thm:simple-power-probabilistic-profile},
$q(H^m)=q(H)$. In particular,
$q_{\cG}=q(H)$. The result now follows from
\Cref{cor:probabilistic-qary-gv}.
\end{proof}

Finite solvable groups give a broader application of the characterization
in \Cref{thm:GV-characterization}. Mann proved that every finite solvable
group $G$ satisfies
$\mathcal M(G)\leq13\rk(G)/4+2$
\cite[proof of Theorem~10, p.~437]{Mann1996}. Lucchini, Menegazzo, and
Morigi obtained the sharp leading constant as follows. Let
$\widehat F_r^{\mathrm{sol}}$ be the free prosolvable group of rank $r$.
For a finitely generated profinite group $\Gamma$, let
\[
m_s(\Gamma)
:=
\#\{M<_{\max}\Gamma:M\text{ is open and }[\Gamma:M]=s\},
\]
and 
$
\mathcal M(\Gamma)
:=
\sup_{\substack{s\geq2\\m_s(\Gamma)>0}}
\frac{\log m_s(\Gamma)}{\log s}.
$
Their theorem states that, for $r\geq3$,
\[
c_3(r-1)+1
\leq \mathcal M(\widehat F_r^{\mathrm{sol}})
\leq
\max\{c_3(r-1)+1,c_2r\},
\]
where
$c_3:=1+\log_9 48+\frac13\log_9 24
=3.243991\ldots$ and
$c_2:=1+\log_4 6+\frac13\log_4 24=3.056\ldots$
\cite[Theorem~4]{LucchiniMenegazzoMorigi2006}.
The upper-bound argument leading to this theorem gives, for every $s\geq2$,
$m_s(\widehat F_r^{\mathrm{sol}})\leq
s^{c_3r+2}$
\cite[proof of Theorem~4]{LucchiniMenegazzoMorigi2006}.

Now let $G$ be a finite solvable group with $\rk(G)=r$. By the universal
property, $G$ is a quotient of $\widehat F_r^{\mathrm{sol}}$. Taking
preimages sends distinct maximal subgroups of $G$ of index $s$ to distinct
maximal open subgroups of $\widehat F_r^{\mathrm{sol}}$ of the same index.
Consequently, $m_s(G)\leq m_s(\widehat F_r^{\mathrm{sol}})$, and hence
\begin{equation}
\label{eq:solvable-maximal-subgroup-bound}
\mathcal M(G)\leq c_3\rk(G)+2.
\end{equation}
Thus every family of nontrivial finite solvable groups has uniform PMSG
and hence is a GV family by \Cref{thm:GV-characterization}. When the
generating ranks tend to infinity,
\eqref{eq:solvable-maximal-subgroup-bound} gives
$\mu_{\cG}\leq c_3$, and \Cref{thm:GV-range} gives the
following explicit rate.

\begin{corollary}[Finite solvable groups]
\label{cor:solvable-groups-gv}
Let \(\cG=(G_i)_{i\geq1}\) be a family of nontrivial finite solvable
groups such that \(\rk_{\cG}=\infty\). Then, for every
\(0<\delta<1-1/q_{\cG}\),
\[
R_{\cG}(\delta)
\geq
\frac{1}{c_3}
\bigl(1-H_{q_{\cG}}(\delta)\bigr),
\]
where, as usual, \(H_{q_{\cG}}(\delta)=\delta\) if
\(q_{\cG}=\infty\).
\end{corollary}

The lower bound in the quoted theorem also shows that the coefficient
$c_3$ is sharp. Indeed, with the analogous definition of $\mu$ for profinite
families using maximal open-subgroup growth, the family
$\widehat\cF^{\mathrm{sol}}=(\widehat F_r^{\mathrm{sol}})_{r\geq3}$
satisfies $\mu_{\widehat\cF^{\mathrm{sol}}}=c_3$.

Using \eqref{eq:solvable-maximal-subgroup-bound} in
\Cref{lem:entropy-sensitive-failure} also shows
that every nontrivial finite solvable group $G$ with $\rk(G)\geq8$ has a generating multiset
$A$ with $|A|\leq7\rk(G)$ and $\delta(A;G)\geq1/10$. Thus,
the constant $85$ in \cite[Theorem~1.13]{ChapmanDinurLubotzky} can be
replaced by $7$.

After showing that many families are GV, we conclude with a family that is
not GV.
\begin{example}
\label{exam:alternating-factorial-power}
Let $G_\ell=A_\ell^{\ell!}$ for $\ell\geq7$. By
\cite[Theorem~1]{Dixon1969} and
\cite[Corollary~7]{Kantor1990ThePO},
we have $\rk(G_\ell)=3$ for all sufficiently large $\ell$.

For each $j\in[\ell!]$, the subgroup
\[
M_{\ell,j}
=
A_\ell^{j-1}\times\operatorname{Stab}_{A_\ell}(1)
\times A_\ell^{\ell!-j}
\]
is maximal in $G_\ell$ and has index $\ell$. Hence $G_\ell$ has at least
$\ell!$ maximal subgroups of index $\ell$.
For every fixed $C>0$,
$
\frac{\ell!}{\ell^{C \rk(G_\ell)}}
=
\frac{\ell!}{\ell^{3C}}
\longrightarrow\infty.
$
Thus $(G_\ell)_{\ell\geq7}$ does not have uniform PMSG and, by
\Cref{thm:GV-characterization}, is not GV.
\end{example}
\printbibliography
\end{document}